\documentstyle{article}
\begin{document}
\bigskip

\begin{center}
{\bf Numerical Approximation of Real
Functions and One Minkowski's Conjecture on Diophintine
Approximations  }

\end{center}
\smallskip
\begin{center}
{\bf Nikolaj M. Glazunov} \end{center}
\smallskip
\begin{center}
{\rm Glushkov Institute of Cybernetics NAS \\
 03187 Ukraine Kiev-187 Glushkov prospekt 40 \\
 Email:} {\it glanm@d105.icyb.kiev.ua }
\end{center} \smallskip

\begin{center} {\bf  Abstract} \end{center}

   In this communication  I consider the applications of
several kinds of approximations of real functions  to the problem
of verified computation  (reliable computing) of the range of
implicitly defined real function
$ x_{n+1} = G(x_{1}, \cdots, x_{n}), $ where dependency
$ F(x_{1}, \cdots, x_{n+1}) = 0 $ is
defined on some compact domain by a sufficiently smooth real
function $ F(x_{1}, \cdots, x_{n+1}) .$
Constructive version of Kolmogorov-Arnold and implicit function
theorems, results about floating-point approximation,
floating-point  approximations which give lower-bound and
upper-bound estimates of some real functions, and approximate
algebraic computation are used for the purpose.  The rigorous
theory  can be build on the base of analysis on manifolds
over floating points domains.
In the text we demonstrate our approach on examples.

\begin{center} {\bf Introduction} \end{center}
In~\cite{AAC} the notion of Approximate Algebraic Computation
 (AAC) is formulated and shortly discussed.
The subject matter of this talk lies in the area between the
geometry of numbers and the analysis of real functions and
approximate algebraic computation of the functions.
More specifically I want to discuss approximate algebraic
computation  aspects of relation
between the geometry of parametric minima of convex and distance
functions and the analysis of functions which determine critical
determinants.  Let
$$ |\alpha x + \beta y|^p + |\gamma x + \delta y|^p \leq c
 |\det(\alpha \delta - \beta \gamma)|^{p/2}, $$
be a diophantine inequality defined for a given real $ p >1 $;
hear $\alpha, \beta, \gamma, \delta$ are real numbers with
$ \alpha \delta - \beta \gamma \neq 0 .$
H. Minkowski in his monograph~\cite{Mi:DA} raise the question
about minimum constant $c$ such that the inequality has integer
solution other than origin.
 This Minkowski's problem can be reformulated as a
conjecture concerning the critical determinant of the region
 $ \mid x \mid^p + \mid y \mid^p  \ \leq 1, \ p > 1.$
 Mentioned mathematical problems are closely connected with
Diophantine Approximation.  \\
Also I want to proposed algorithms for approximate algebraic
computation. By the computation of the algorithms A. Malishev
and I have investigated  the Minkowski's conjecture and proposed
the strengthen Minkowski's analytic (MAS) conjecture.
For verified computation (reliable computing) we used interval
analyses. Methods and algorithms for interval evaluation of
explicitly and implicitly defined real functions are used.
In  the paper we briefly consider following topics:                \\
critical determinant of a body and diophantine approximation;  \\
the problem;          \\
Minkowski's analytic conjecture;    \\
Kolmogorov-Arnold's theorem; \\
computer experiments and
strengthen Minkowski's analytic conjecture;   \\
interval-analytic methods; \\
algorithms for approximate algebraic computation
and implementation.\\

\begin{center} {\bf 2 Critical Determinant of a Body
and Diophantine Approximation} \end{center}

  Critical determinant is one of the main notion of the Geometry
of Numbers. Recall the definitions~\cite{C:GN}.
Let $\cal D $ be a set and $\Lambda $ be a  lattice with base
$ \{a_1, \ldots ,a_n \}$ in ${\bf R}^n.$ A lattice $\Lambda $
is {\em admissible} for body $\cal D $ ($ {\cal D}-{\em admissible}$)
if ${\cal D} \bigcap \Lambda = \emptyset $ or $0.$
Let $ d(\Lambda) $ be the determinant of $\Lambda.$ The infimum
$\Delta(\cal D) $ of determinants of all lattices admissible for
$\cal D $ is called {\em the critical determinant} of $\cal D; $
if there is no $\cal D-$admissible lattices then puts
$\Delta(\cal D) = \infty. $ A lattice $\Lambda $ is {\em critical}
if $ d(\Lambda) = \Delta(\cal D).$ \\
Critical determinant is closely connected with diophantine
approximation, solving inequalities $ F(x_1, \ldots ,x_n) < c$ in
integer numbers $ x_1, \ldots ,x_n $ (with some restrictions,
for instance, $ x = (x_1, \ldots ,x_n) \neq 0 $). Usually in the
geometry of numbers the function $F(x)$ is a distance function.
 A real function $F(x)$ defined on
${\bf R}^n$ is {\em distance function} if \\
  (i) $F(x) \geq 0, x \in {\bf R}^n, F(0) = 0;$ \\
 (ii) $F(x)$ is continuous;   \\
(iii) $F(x)$ is homogenous: $F(\lambda x) = \lambda F(x), \lambda
\in {\bf R}$. \\
The problem of solving of diophantine inequality $F(x) < c $,
with a distance function $F$ are investigated. \\
Let $\overline M$ be the closure of a set $M$ and $\#P$ be the number
of elements of a finite set $P$. An open set $ S \subset {\bf R}^n $
is a {\em star body} if $S$ includes the origin of ${\bf R}^n$ and for
any ray $r$ beginning in the origin
$ \#(r \cap (\overline M \setminus M)) \leq 1$. If $F(x)$ is a
distance function then the set
$$ M_F = \{ x: F(x) < 1\} $$
is a star body.             \\
One of the main particular case of a distance function is the case
of convex symmetrical function $F(x)$ which with conditions (i) - (iii)
satisfies the additional conditions \\
(iv)$ F(x + y) \leq F(x) + F(y) ;$  \\
(v) $F(-x) = F(x) .$

\begin{center} {\bf 3 The Problem} \end{center}

In considering the question of the minimum value taken by the
expression $ |x|^p + |y|^p $, with $ p \geq 1 $, at points, other
that the origin, of a lattice $\Lambda$ of determinant
$ d(\Lambda)$, Minkowski~\cite{Mi:DA} shows that the problem
of determining the maximum value of the minimum for different
lattices may be reduced to that of finding the minimum possible
area of a parallelogram with one vertex at the origin and the
three remaining vertices on the curve $ |x|^p + |y|^p = 1 $.
The problem with $ p = 1, 2$ and $\infty$ is trivial: in these
cases the minimum areas are $1/2, \: \sqrt{3}/2 $ and $1$
respectively.
 Let $ D_p \subset {\bf R}^2 = (x,y), \ p > 1 $ be the 2-dimension
region:

$$ |x|^p + |y|^p < 1 . $$
Let $\Delta(D_p) $ be the critical determinant of the region.
 Using analytic parameterization Cohn~\cite{Co:MC} gives analytic
formulation of Minkowski's conjecture.
 Let
$$ \Delta(p,\sigma) = (\tau + \sigma)(1 + \tau^{p})^{-\frac{1}{p}}
  (1 + \sigma^p)^{-\frac{1}{p}}, \; \; \;  \; (1) $$

be the function defined in the domain
 $$ D_{p}: \; \infty > p > 1, \; 1 \leq \sigma \leq \sigma_{p} =
 (2^p - 1)^{\frac{1}{p}}, $$

of the $ \{p,\sigma\} $ plane, where $\sigma$ is some real parameter;
$\;$ here $ \tau = \tau(p,\sigma) $ is the function uniquely
determined by the conditions
$$ A^{p} + B^{p} = 1, \; 0 \leq \tau \leq \tau_{p}, $$
where
$$ A = A(p,\sigma) = (1 + \tau^{p})^{-\frac{1}{p}} -
(1 + \sigma^p)^{-\frac{1}{p}}                  $$
$$ B = B(p,\sigma) = \sigma(1 + \sigma^p)^{-\frac{1}{p}}       +
\tau(1 + \tau^{p})^{-\frac{1}{p}}, $$
$\tau_{p}$ is defined by the equation
$$ 2(1 - \tau_{p})^{p} = 1 + \tau_{p}^{p}, \; 0 \leq \tau_{p}
\leq 1. $$
In this case needs to extend the notion of parameter variety to
parameter manifold. The function $ \Delta(p,\sigma) $ in region
$ D_{p} $  determines the parameter manifold.

{\bf Minkowski's analytic conjecture:} \\

{ \it For any real $p$ and $\tau$ with conditions }
$ p > 1, \ p \ne 2, \ 0 < \tau < \tau_p $ \\
$$ \Delta(p,\sigma) > min(\Delta(p,1),\Delta(p,\sigma_p)).$$ \\

 For investigation of properties of function $  \Delta(p,\sigma) $
which are need for proof of Minkowski's conjecture~\cite{Mi:DA,Co:MC}
we considered the value of $\Delta = \Delta(p,\sigma) $ and its
derivatives $  \Delta_{\sigma}^{'} \; , \,
\Delta_{\sigma^{2}}^{''} \; , \; \Delta_{p}^{'} \;, \;
 \Delta_{\sigma p}^{''} \; , \; \Delta_{\sigma^{2}p}^{'''} \;$
on some subdomains of the domain  $ D_{p} $~\cite{Gl:CA,Gl:VV,
Gl:PT,GGM:PM}. The analytical computation of the derivatives is a
problem of computer algebra.

\begin{center} {\bf 4 The Theorem of Kolmogorov-Arnold}
\end{center}

  Let $ F(x_{1}, \cdots, x_{n+1}) = 0, n > 1 $ be a sufficiently
smooth real function of two variables. The rough form of the
theorem of Kolmogorov-Arnold states that: \\

{\bf Theorem} (Kolmogorov, Arnold)          \\
{\it Any sufficiently smooth real function can be represented as a
superposition of functions of two variables.} \\

  Let us demonstrate the constructive version of the theorem
on examples of function $\Delta = \Delta(p,\sigma,\tau) $ and it's
derivatives. At first expressing $  \Delta_{\sigma}^{'} \; , \;
\Delta_{\sigma^{2}}^{''} \; , \; \Delta_{p}^{'} \; , \;
 \Delta_{\sigma p}^{''} \; , \; \Delta_{\sigma^{2}p}^{'''} \; $
in terms of a sum of derivatives of "atoms" $ s_{i} = \sigma^{p-i},
\; t_{i} = \tau^{p-i}, \; a_{i} = (1 + \sigma^{p})^{-i-\frac{1}{p}},
\; b_{i} = (1 + \tau^{p})^{-i-\frac{1}{p}}, \; A = b_{0} - a_{0},
\; B = \tau b_{0} + \sigma a_{0}, \; \alpha_{i} = A^{p-i}, \;
\beta_{i} = B^{p-i} \; ( i = 0, 1, 2, \ldots).$
  Then by the implicit function theorem computing
$ \tau = \tau(p,\sigma) $ by means of the following iteration
 process:

$$ {\tau}_{i + 1} = (1 + {\tau}_{i}^p)^{\frac{1}{p}}
((1 - ((1 + {\tau}_{i}^p)^{-\frac{1}{p}} - (1 + {\sigma}^{p})^
{-\frac{1}{ p}})^{ p})^{\frac{1}{ p}} -
{ \sigma}(1 + {\sigma}^{p})^{-\frac{1}
{ p}}), $$

  For approximate computation of the expression for $\tau_p $
we apply the following iteration: \\

$$ {(\tau_p)}_{i + 1} = 1 - (2^{-\frac{1}{p}})(1 + {(\tau_p)}_{i}^p)
^{\frac{1}{p}}, \; p > 1, \; {(\tau_p)}_{0} \in [0,0.36].  $$

So we really have represented the function $\Delta$  as the
function $\Delta(p,\sigma)$ of two variables. The same fact
is true for it's derivatives
  Now we can compute expressions for $\Delta, \Delta_{\sigma}^{'}
 \; , \; \Delta_{\sigma^{2}}^{''} \; , \; \Delta_{p}^{'} \; , \;
 \Delta_{\sigma p}^{''} \; , \; \Delta_{\sigma^{2}p}^{'''} \; $
by means of approximate algebraic computations.

\begin{center} {\bf 5 Strengthen Minkowski's analytic conjecture }
\end{center}

Based on some theoretical evidences and results of mentioned
computation A.V. Malishev and author proposed \\

{\em Strengthen Minkowski's analytic (MAS) conjecture:} \\

{\it For given $ p > 1 $ and increasing $\sigma$ from $0$ to
$\sigma_p$ the function $\Delta(p,\sigma) $ \\
1) increase strictly monotonous if $ 1 < p < 2$ and
$ p \geq p^{(1)} $, \\
2) decrease strictly monotonous if $ 2 \leq p \leq p^{(2)} $, \\
3) has a unique maximum on the segment $ (1,\sigma_p) $; until the
maximum $\Delta(p,\sigma) $ increase strictly monotonous and
then decrease strictly monotonous if $ p^{(2)} < p < p^{(1)} $; \\
4) constant, if $ p = 2$, \\
here $ p^{(1)} > 2 $ is a root of equation
$ \Delta_{\sigma^{2}}^{''}|_{\sigma = \sigma_p } = 0 $; $
p^{(2)} > 2 $ is a root of equation
$ \Delta_{\sigma^{2}}^{''}|_{\sigma = 1} = 0 $ }. \\
It is seems that conjecture (MAS) did not proven  for any
parameter $p$ except trivial $ p = 2 $.  \\

\begin{center} {\bf 6 Interval-analytic methods }
\end{center}

Let $ {\bf X} = ({\bf x}_{1}, \cdots,{\bf x}_{n}) =
([{\underline x}_{1}, {\overline x}_{1}], \cdots,
[{\underline x}_{n}, {\overline x}_{n}] $ be the n-dimensional
real interval vector with
$ {\underline x}_{i} \leq x_{i} \leq {\overline x}_{i} $
 ("rectangle" or "box"). The {\it interval evaluation} of a
function $ G(x_{1}, \cdots, x_{n}) $ on an interval ${\bf X} $
is the interval $[{\underline G}, {\overline G}] $ such that
for any $ x \in {\bf X}, \; G(x) \in [{\underline G}, {\overline G}]. $
The interval evaluation is called {\it optimal} if
$ {\underline G} = \min G,$ and $ {\overline G} = \max G $
on the interval {\bf X}.   \\
 In the communication I consider the case $ n = 2 \; $\footnote{By
the result of A. Kolmogorov and V. Arnold, any sufficiently
smooth real function can be represented as a superposition
of functions of two variables.} It is sufficient for Minkowski's
conjecture. For the purpose we used modified variant of the
method~\cite{Ma:AC} which we called \\
 {\bf Malyshev's Method:} \\
Let $ D $ be a subdomain of  $ D_{p}. $ Under evaluation in
$ D $ a mentioned function the domain
is covered by rectangles of the form
$$ [{\underline p}, {\overline p}; \;
{\underline \sigma}, {\overline \sigma}]. $$
In the case of formula (1) expressing $  \Delta_{\sigma}^{'} \; , \;
\Delta_{\sigma^{2}}^{''} \; , \; \Delta_{p}^{'} \; , \;
 \Delta_{\sigma p}^{''} \; , \; \Delta_{\sigma^{2}p}^{'''} \; $
in terms of a sum of derivatives of "atoms" $ s_{i} = \sigma^{p-i},
\; t_{i} = \tau^{p-i}, \; a_{i} = (1 + \sigma^{p})^{-i-\frac{1}{p}},
\; b_{i} = (1 + \tau^{p})^{-i-\frac{1}{p}}, \; A = b_{0} - a_{0},
\; B = \tau b_{0} + \sigma a_{0}, \; \alpha_{i} = A^{p-i}, \;
\beta_{i} = B^{p-i} \; ( i = 0, 1, 2, \ldots).$ one applies the
rational interval evaluation to construct formulas
for lower bounds and upper bounds of the functions, which in the
end can be expressed in terms of $ {\underline p}, \; {\overline p}, \;
{\underline \sigma}, \; {\overline \sigma}, \; {\underline \tau},
\; {\overline \tau}, \; ; \; $ here the bounds $ \; {\underline \tau},
 \; {\overline \tau}, \; $ are obtained with the help of the iteration
process:

$$ {\underline t}_{i + 1} = (1 + {\underline t}_{i}^{\overline p})^
{\frac{1}{\overline p}}((1 - ((1 + {\underline t}_{i}^{\overline p})^
{-\frac{1}{\overline p}} - (1 + {\overline \sigma}^{\underline p})^
{-\frac{1}{\underline p}})^{\underline p})^{\frac{1}{\underline p}} -
{\overline \sigma}(1 + {\overline \sigma}^{\underline p})^{-\frac{1}
{\underline p}}), $$

$$ {\overline t}_{i + 1} = (1 + {\overline t}_{i}^{\underline p})^
{\frac{1}{\underline p}}((1 - ((1 + {\overline t}_{i}^{\underline p})^
{-\frac{1}{\underline p}} - (1 + {\underline \sigma}^{\overline p})^
{-\frac{1}{\overline p}})^{\overline p})^{\frac{1}{\overline p}} -
{\underline \sigma}(1 + {\underline \sigma}^{\overline p})^{-\frac{1}
{\overline p}}). $$
           $$  \; i = 0,1,\cdots  $$
As interval computation is the enclosure method, we have to put:

$$ [{\underline \tau}, \; {\overline \tau}] =
 [{\underline t}_N, \; {\overline t}_N] \bigcap
 [{\underline \tau}_{0}, \; {\overline \tau}_{0}] \; .$$
$ N $ is
computed on the last step of the iteration. \\
For initial values we may take $: \; [{\underline t}_{0}, \;
{\overline t}_{0}] =
[{\underline \tau}_{0}, \; {\overline \tau}_{0}] = [0,\; 0.36]. $

 \begin{center} {\bf 7 Algorithms for approximate algebraic
 computation } \end{center}

 Below each of 6 first algorithms has 4 different forms: \\
 (i) approximate algebraic computation of the given
 expression in a given point at the floating point representation; \\
 (ii) approximate algebraic computation of a lower-bound estimate
 of the given expression in a given point at the floating point
 representation; \\
 (iii) approximate algebraic computation of a upper-bound estimate
 of the given expression in a given point at the floating point
 representation; \\
 (iv) approximate algebraic computation of an interval
 evaluation  of the given expression over given intervals at the
 floating point representation; \\
  Here we give names, input and output of  algorithms
  for interval evaluation only. But all algorithms  for Minkowski
conjecture are  implemented and tested. \\

 {\bf Algorithm} {\em TPV} \\
{\bf Input:} An implicitly defined function $\tau_p$ from
Paragraphs 3 and 6. \\
$ [{\underline p}, {\overline p}; \; {\underline \sigma},
{\overline \sigma}].$  \\
{\bf Method:} Shortly described in Paragraphs 3 and 6. \\
{\bf Output:} The interval evaluation of $\tau_p$.    \\

{\bf Algorithm} {\em TAUV} \\
{\bf Input:} Implicitly defined function $\tau$ from
Paragraphs 3 and 6. \\
$ [{\underline p}, {\overline p}; \; {\underline \sigma},
{\overline \sigma}].$  \\
{\bf Method:} Shortly described in Paragraphs 3 and 6. \\
{\bf Output:} The interval evaluation of $\tau$.   \\

{\bf Algorithm} {\em L0V} \\
{\bf Input:}  Function $l^0 =  \Delta(p,\sigma) - \Delta_p^{(0)}$
from Paragraphs 3 and 6. \\
$ [{\underline p}, {\overline p}; \; {\underline \sigma},
{\overline \sigma}].$  \\
{\bf Method:} Shortly described in Paragraphs 3 and 6. \\
{\bf Output:} The interval evaluation of $l^0.$    \\

{\bf Algorithm} {\em L1V} \\
{\bf Input:}  Function $l^1 = \Delta(p,\sigma) - \Delta_p^{(1)}$
from Paragraphs 3 and 6. \\
$ [{\underline p}, {\overline p}; \; {\underline \sigma},
{\overline \sigma}].$  \\
{\bf Method:} Shortly described in Paragraphs 3 and 6. \\
{\bf Output:} The interval evaluation of $l^1.$    \\

{\bf Algorithm} {\em GV} \\
{\bf Input:}  A function $g(p,\sigma)$  which has the same sign
as function $\Delta_{\sigma}^{'}$.\\
$ [{\underline p}, {\overline p}; \; {\underline \sigma},
{\overline \sigma}].$  \\
{\bf Method:} Shortly described in Paragraphs 3 and 6. \\
{\bf Output:} The interval evaluation of $g(p,\sigma).$    \\

{\bf Algorithm} {\em HV} \\
{\bf Input:}  A function $h(p,\sigma)$  which is the partial
derivative by $\sigma$ the function $g(p,\sigma)$.\\
$ [{\underline p}, {\overline p}; \; {\underline \sigma},
{\overline \sigma}].$  \\
{\bf Method:} Shortly described in Paragraphs 3 and 6. \\
{\bf Output:} The interval evaluation of $h(p,\sigma).$    \\

Next two algorithms are described in~\cite{Gl:IA}.

{\bf Algorithm} {\em MonotoneFunction} \\
{\bf Input:} A real function $F(x,y)$ monotonous by $x$ and
by $y.$ \\
Interval $[{\underline x},{\overline x}; {\underline y},
{\overline y}].$ \\
{\bf Output:} The interval evaluation of $F.$    \\

{\bf Algorithm} {\em RationalFunction} \\
{\bf Input:} A rational function $R(x,y).$
Interval $[{\underline x},{\overline x}; {\underline y},
{\overline y}].$ \\
{\bf Output:} The interval evaluation of $R.$    \\

\begin{center} {\bf Implementation} \end{center}

 Algorithms are implemented on PL/1(O) and (partly)
on C and C++.
There is also an implementation of some of the algorithms on
Reduce3.3-3.X in T. Sasaki's package~\cite{Sa}
of arbitrary precision rational arithmetic.

\begin{center} {\bf Acknowledgments} \end{center}

Author thanks to Prof. Matu-Tarow Noda for his proposition
to present the text as an electronic publication in the
frame of the IMACS-ACA session "Approximate Algebraic
Computation".

\end{document}